\documentstyle[amssymb,amsmath,eucal]{article}

\begin{document}

\def\D{{\rm D}}

\def\R {{\Bbb R }}
\def\C {{\Bbb C }}
\def\KK{{\Bbb K}}
\def\HH{{\Bbb H}}

\def\U{{\rm U}}

\def\const{{\rm const}}
\def\B{{\rm B}}

\def\OO{{\rm O}}
\def\SO{{\rm SO}}
\def\GL{{\rm GL}}
\def\SL{{\rm SL}}
\def\SU{{\rm SU}}
\def\Sp{{\rm Sp}}
\def\SOS{\SO^*}

\def\Mat{{\rm Mat}}

\def\Subsection{}
\def\Subsections{}

\newcommand{\rk}{\mathop{\rm rk}\nolimits}

\def\ov{\overline}
\def\phi{\varphi}
\def\epsilon{\varepsilon}
\def\kappa{\varkappa}
\def\le{\leqslant}
\def\ge{\geqslant}

\newcommand{\sh}{\mathop{\rm sh}\nolimits}
\newcommand{\ch}{\mathop{\rm ch}\nolimits}

\renewcommand{\Re}{\mathop{\rm Re}\nolimits}
\renewcommand{\Im}{\mathop{\rm Im}\nolimits}
\newcommand{\tr}{\mathop{\rm tr}\nolimits}

\def\matritsa{
\begin{pmatrix}a&b\\ \ov b&\ov a\end{pmatrix}}

\begin{center}

{\large\bf
Action of overalgebra in Plancherel
decomposition and shift operators
 in imaginary direction}

 \smallskip

{\sc Yurii A. Neretin }\footnote{Supported by
the grant NWO 047--008--009}

\end{center}

In last 50 years,
there were solved
many problems of
explicit spectral decompositions
for restrictions of unitary representations
to subgroups (see the bibliography
in \cite{MolVINITI}).

This class of problems
 also contains  the following problems,
which are formulated
in other terms.

1. {\it Decomposition of the tensor product
of representations
$\rho_1$, $\rho_2$
of a group $G$}.
Indeed, this is exactly the problem
of restriction of representations of $G\times G$
to the diagonal subgroup $G\subset G\times G$.

\smallskip

2.
{\it Decomposition of $L^2$
on a pseudo-Riemannian symmetric space
 $G/H$}.
As was observed\footnote{This phenomenon
could be easily observed
from the Makarevich paper \cite{Mak}
published in 1973,
in a strange way it remained
 nonformulated
for a long time}
 in \cite{pse},
 for each classical pseudo-Riemannian
symmetric space  $G/H$,
there exists a canonical classical group
 $\widetilde G\supset G$
and a representation
$\rho$ of the group $\widetilde G$
of degenerated principal series
satisfying one of two following properties
(usually the first variant is realized):

--- the restriction of
$\rho$ to $G$ is  $L^2(G/H)$

--- the restriction of
$\rho$ to $G$ is
the  direct sum
of the spaces
$L^2(G/H_j)$, there $G/H_j$
is a finite collection of   symmetric spaces,
and $G/H$ is one of the spaces $G/H_j$.

Hence the decomposition of $L^2$
on a classical Riemannian symmetric space
can be considered as a restriction problem.

Description of the spectral type
(without explicit Plancherel formula)
for all the pseudo-Riemannian
symmetric spaces was recently obtained
in the works of van den Ban,
Schlichtkrull,  Delorme,
and Oshima
(the proof is contained in the union of
a large collection of papers,
for references, see \cite{Ban}, \cite{Del}).
It seems that for classical
symmetric spaces
the problem of evaluation
of the Plancherel measure is near the final solution.

For some cases, the explicit Plancherel
measure is known; in particular,
for
 $L^2$ on semisimple groups (\cite{GN}, \cite{HC}),
on Riemannian symmetric spaces
   (\cite{GN},
\cite{GK1}, \cite{GK2}, see also \cite{Hel}),
  on rank 1 spaces  (\cite{Molhyper}, \cite{Molrank})
and on the spaces
 $G_\C/G_\R$, where $G_\C$ is a complex group,
and $G_\R$
is its real form  (\cite{Har1}, \cite{Har2}).

\smallskip

3. {\it Berezin kernel representations}
 (deformations of $L^2$
on Riemannian noncompact symmetric space $G/K$)
also can be obtained by the restriction
from some overgroup $G^\circ\supset G$,
see \cite{Nerpl}, \cite{Nerhyper}
and references in \cite{Nerhyper}.

\smallskip

Usually the problem of the noncommutative
harmonic analysis is formulated as
the problem about the spectrum of a representation
or as the more complicated problem of an explicit
decomposition of a representation
into a direct integral of irreducible
representations;
the last question
includes
 the explicit evaluation of the spectral
measure (the so-called
Plancherel measure).
Recent works of the author
 \cite{Nerindex}, \cite{Nerhyper}
contain an attempt
of an investigation of "the analysis
after the Plancherel formula".

This work is a continuation of
 \cite{Nerindex}.
 Here we are trying to understand
  the answer to following
 question.

\smallskip

{\sc Question.} Assume that we know the explicit
 Plancherel formula for the restriction
 of a unitary representation $\rho$
 of a group  $G$
to a subgroup $H$.
Is it possible to write the action of Lie algebra
of
 $G$ in the direct integral of  representations
 of
$H$?

\smallskip

We obtain the positive solution
of this problem for one of the simplest possible
examples, precisely, for the tensor product
of a representation of $\SL(2,\R)$
with a highest weight and the conjugate representation
with a lowest weight.
This tensor product and its decomposition
were widely discussed in the literature
on representation theory and special
functions in last 40 years (some references: \cite{Puk},
\cite{dB}, \cite{VGG}, \cite{tensor},
\cite{vD}). Nevertheless,
the formula for the action of $sl_2\oplus sl_2$
in tensor product were not appear.
The reason is the unusual for representation
theory form of these Lie algebra operators.

It turns out  that the operators
(12)--(14)
of Lie algebra $sl_2\oplus sl_2$
 are the second order differential
operators with respect to one variable
and the second order difference operators with respect
to the another variable; moreover,
it turns out that the difference operators are
defined in the terms of a shift
in the imaginary direction, i.e.,
in our formula,
there appear the operators
of the form
$$ T f(x)=f(x+i) \qquad \qquad \text{for $f\in L^2(\R)$}$$
($i^2=-1$).
There is no self-contradiction in this expression,
the shift operators are well defined
on functions admitting the holomorphic continuation
to an appropriate strip.
Operators $f\mapsto xf$ and $f\mapsto \frac d{dx}$
also are not defined on the whole $L^2(\R)$.

First, these Lie algebra  operators were obtained by
author using theorems on the
operational calculus for the index hypergeometric
transform (17) from \cite{Nerindex}.
But the final formulas are elementary
and admit a direct  verification
by elementary tools.

Since  the approximate
 structure of formulas for action
 of overalgebra now became more clear,
 it is natural to formulate the general
 problem given above.


\medskip

This work was prepaired
during my visit to the Erwin Schr\"odinger
Institute (Vienna).
I thank the administrators of the Institute
for their hospitality.

\smallskip

{\bf 1. Group $\SL(2,\R)$.}
We realize the group $\SL(2,\R)$
as the group of complex
 $2\times 2$ matrices having the form
$$
\begin{pmatrix}a&b\\
\ov b&\ov a\end{pmatrix}, \qquad \text{where}\quad
  |a|^2-|b|^2=1
.$$

By $\D$
we denote the disk $|z|<1$
on the complex plane
$\C$, by $S^1$
we denote the circle $|z|=1$;
we  represent the points
of the circle in the form  $z=e^{i\phi}$.

The group $\SL(2,\R)$
acts on the disk by
the M\"obius
transformations
\begin{equation}
\begin{pmatrix}a&b\\ \ov b&\ov a\end{pmatrix}: \quad
z\mapsto \frac{az+b}{\ov b z+\ov a}
.\end{equation}
If $z\in S^1$, then its image
under (1) also
is contained in $S^1$.

\smallskip

{\bf 2. Highest weight representations
of $\SL(2,\R)$.} Fix
$\alpha>1$.
Consider the space $H_\alpha$
of holomorphic functions  $f$
 on the circle $\D$
 satisfying the condition
\begin{equation}
\int_\D|f(z)|^2 (1-|z|^2)^{\alpha-2}\{dz\}<\infty
,\end{equation}
where $\{dz\}$ denotes the Lebesgue measure
on $\D$. Define
the inner product in the space $H_\alpha$
by
\begin{equation}
\langle f, g\rangle_\alpha=\frac{\alpha-1}\pi
 \int_\D f(z) \ov{g(z)} (1-|z|^2)^{\alpha-2}\{dz\}
;\end{equation}
 the pre-integral factor
 provides the identity
 $\langle 1,1\rangle_\alpha=1$.
 In terms of the Taylor coefficients,
 the inner product is given by
$$
\langle \sum c_k z^k, \sum c_k'z^k\rangle_\alpha=
\sum_{k\ge 0} c_k \ov c_k'\frac{k!}
{\alpha(\alpha+1)\dots (\alpha+k-1)}
,$$
and condition (2) for $f(z)=\sum c_k z^k$
has the form
\begin{equation}
\sum |c_k|^2 k^{-(\alpha-1)}<\infty
.\end{equation}

It can easily be checked that the space $H_\alpha$
 is complete with respect to the  norm
 defined by this inner product,
 i.e., $H_\alpha$ is a Hilbert space.

The group $\SL(2,\R)$
acts in the space $H_\alpha$
by the unitary operators
\begin{equation}
T_\alpha\matritsa(z)=f\Bigl(\frac{az+b}
{\ov b z+\ov a}\Bigr)
  (\ov b z+\ov a)^{-\alpha}
.\end{equation}

First, consider the case, when $\alpha$
is an integer.
The factor $(\ov b z+\ov a)^{-\alpha}$
 is a degree of the derivative of the function (1),
this implies that $T_\alpha$
is a representation:
\begin{equation}
T_\alpha(g_1)T_\alpha(g_2)=T_\alpha(g_1g_2)
.\end{equation}
The operators $T_\alpha(g)$ are unitary, i.e.,
$$
\langle T_\alpha(g) f_1, T_\alpha(g) f_2\rangle_\alpha
 = \langle f_1, f_2\rangle_\alpha,
$$
 this can be easily checked by a change of
the variable.

For a noninteger   $\alpha$,
$T_\alpha(g)$ also is an unitary representation;
we only must explain the
meaning of the expression
\begin{equation}
(\ov b z+\ov a)^{-\alpha}=
(1+b\ov a^{\,-1} z)^{-\alpha} a^{-\alpha}=
(1+b\ov a^{\,-1} z)^{-\alpha}
 e^{-\alpha(\ln a+ 2\pi i k)}
.\end{equation}
Obviously, $|b\ov a^{\,-1}|<1$.
Hence the function
$$(1+b\ov a^{-1} z)^{-\alpha}:=
1+\frac{\alpha}{1!}b\ov a^{-1} z+
\frac{\alpha(\alpha-1)}{2!}
(b\ov a^{-1})^2 z^2+\dots$$
is well defined.

Therefore the operator (5)
is defined up to the factor $e^{-2\pi i k\alpha}$,
 the absolute value of this factor is 1.
Hence  equality (6)
is replaced by
$$
T_\alpha(g_1)T_\alpha(g_2)=
\theta\cdot T_\alpha(g_1g_2)
,$$
where $|\theta|=1$.
Thus,    $T_\alpha(g)$
is an unitary projective
representation of the group
$\SL(2,\R)$.

\smallskip

{\sc Remark.}
Clearly, $T_\alpha(g)$
can be also considered as a linear representation
of the universal covering
group
 of $\SL(2,\R)$.

\smallskip

{\bf 3. Tensor product.}
By $\ov  T_\alpha(g)  $
we denote the representation
complex conjugate to $ T_\alpha(g)$,
it acts in the space
 $\ov H_\alpha$ of
 antiholomorphic functions
 in the circle $|u|<1$
 by the formula
$$
T_\alpha\matritsa f(\ov u)=
f\Bigl(\frac{\ov  a\,\ov u+\ov b}{ b \ov u+ a}\Bigr)
  ( b \ov u+ a)^{-\alpha}.
$$
The inner product in $\ov H_\alpha$
is given by (3).

Consider the space $ H_\alpha\otimes\ov  H_\alpha$.
It consists of functions
 $f(z,\ov u)$ on the bidisk $\D\times\D$,
holomorphic with respect to
the variable $z$ and antiholomorphic
in
$u$; the inner product in $ H_\alpha\otimes\ov  H_\alpha$
is given by
$$
\langle f_1,f_2\rangle=\frac{(\alpha-1)^2}{\pi^2}
\iint_{\D \times \D} f_1(z,\ov u)\ov{ f_2(z,\ov u)}
(1-z\ov z)^{\alpha-2} (1-u\ov u)^{\alpha-2} \{dz\}\, \{du\}
.$$
The group $\SL_2(\R)\times\SL_2(\R)$
acts in $H_\alpha\otimes\ov H_\alpha$
by the operators
$ T_\alpha(g_1)\otimes \ov T_\alpha(g_2)$
given
by
\begin{multline*}
T_\alpha
\begin{pmatrix} a_1&b_1\\\ov b_1&\ov a_1
\end{pmatrix}\otimes
 \ov T_\alpha
\begin{pmatrix} a_2&b_2\\\ov b_2&\ov a_2
\end{pmatrix}
 f(z,u)=\\=
f\Bigl(\frac{a_1z+b_1}{\ov b_1 z+\ov a_1}\, ,\,
\frac{\ov  a_2\,\ov u+ \ov b_2}{ b_2 \ov u+ a_2}\Bigr)
(\ov b_1 z+\ov a_1)^{-\alpha}
     ( b_2 \ov u+ a_2)^{-\alpha}
.\end{multline*}
We restrict this representation
to the diagonal
subgroup $\SL_2(\R)\subset\SL_2(\R)\times\SL_2(\R)$,
this corresponds to the substitution
$a_1=a_2=a$, $b_1=b_2=b$ to the last formula.
This  representation
of  $\SL_2(\R)$ is a linear
(nonprojective)
 representation.
Indeed,
$$
(\ov b z+\ov a)^{-\alpha}
     ( b \ov u+ a)^{-\alpha}=
 (1+\ov a^{-1}\ov b z)^{-\alpha}
  (1+ a^{-1}b \ov u)^{-\alpha}
(a\ov a)^{-\alpha}
,$$
and this expression
is a single-valued function
for $|z|<1$, $|u|<1$.

\smallskip

{\bf 4. Principal series of representations.}
 Fix $s\in \R$.
Consider the representation $\rho_s$
of the group $\SL(2,\R)$ in $L^2$
on the circle $S^1$
given by
$$
\rho_s\matritsa f(e^{i\phi})=
f\Bigl(\frac{a e^{i\phi}+b}
{\ov b e^{i\phi}+\ov a}\Bigr)
|\ov b e^{i\phi}+\ov a|^{-1-2is}
.$$
These unitary representations are
the so-called representations
of the principal series.
 Recall, that the representation $T_s$
 is equivalent  to $T_{-s}$
(see  \cite{Bar}).

\smallskip

{\bf 5. Spectral decomposition.} Let $\alpha$
be the same as in
pp.2--3.
Consider the kernel
$$
K_\alpha(\phi, s; z,u):=
\frac{(1-\bar z e^{i\phi})^{-1/2-is}
                          (1- u e^{-i\phi})^{-1/2-is} }
    {(1-\bar z u)^{\alpha-1/2-is}}
,$$
where $|z|<1, |u|<1$, $\phi\in[0,2\pi]$, $s\in\R$.
Consider the integral operator $J_\alpha$,
 that takes
a function $f\in H_\alpha\otimes\ov H_\alpha$
to the function
$F(\phi,s)$
given by
\begin{multline}
F(\phi,s)=
\iint\limits_{\D\times D}
 K_\alpha(\phi, s; z,u)  f(z,\ov u)
  (1-z\ov z)^{\alpha-2}(1-u\ov u)^{\alpha-2} \{dz\}\{du\}
.
\end{multline}

Consider the action of the group $\SL(2,\R)$
 in space of functions in the variables
$(\phi,s)$ defined by
$$
R\matritsa F(e^{i\phi},s)=
F\Bigl(\frac{a e^{i\phi}+b}
{\ov b e^{i\phi}+\ov a},\, s\Bigr)
|\ov b e^{i\phi}+\ov a|^{-1-2is}
.$$
Note, that for a fixed $s\in \R$
the function $F(e^{i\phi},s)$
(as a function in $\phi$)
transforms by the formula
for principal series representations.

Simple calculation show,
that the operator $J_\alpha$
intertwines the representations
$T_\alpha(g)\otimes \ov T_\alpha(g)$ and $R$
of $\SL(2,\R)$:
\begin{equation}
 J_\alpha \cdot \bigl(T_\alpha(g)\otimes
 \ov T_\alpha(g)\bigr)=
   R(g)\cdot J_\alpha
\end{equation}

{\sc Theorem 1.} {\it Operator $J_\alpha$
is a unitary
operator from $H_\alpha\otimes \ov H_\alpha$
to the space
$L^2$ on $\phi\in[0,2\pi]$, $s\ge 0$
with respect to the measure}
\begin{multline}
\frac
{|\Gamma(\alpha-1/2+is)|^2}
{\Gamma(\alpha)^2}\,\, s\,\frac{\sh (\pi s)}
 {\ch (\pi s)}
\,ds\,d\phi=\Bigl|
\frac{\Gamma(\alpha-1/2+is)\Gamma(1/2+is)}
{\Gamma(\alpha)\Gamma(is)}
\Bigr|^2
\,\,
\,ds\,d\phi
.\end{multline}

Thus $T_\alpha\otimes \ov T_\alpha$
is a multiplicity free integral
over principal series representations
(this fact was obtained by Pukanszky \cite{Puk}).
Various ways for obtaining the
Plancherel measure (9)
are contained in \cite{VGG}, \cite{tensor},
 \cite{vD}, \cite{Nerindex}.

\smallskip

{\bf 6. Holomorphic
 continuation of $J_\alpha f(\phi,s)$.}
Denote by $W$ the space
of all $f\in H_\alpha\otimes \ov H_\alpha$
admitting analytic  continuation
to a bidisk
 $|z|<1+\delta$, $|u|<1+\delta$.

For $f\in W$, consider its Taylor series
$$f(z,\ov u)=\sum_{k,l} c_{kl} z^k \ov u^l.$$
Obviously, its coefficients $c_{kl}$
exponentially decrease
for $k+l\to\infty$.

\smallskip

{\sc Lemma.}
{\it Fix  $\phi$.
For $f\in W$, the function $J_\alpha f(\phi,s)$
can be extended holomorphically
to the whole complex plane $s\in\C$}.

\smallskip

 {\sc Proof.}
 Fix
$s\in \C$, $\phi\in[0,2\pi]$.
The function $K_\alpha(\phi,s;  z,u)$
as a function in the variables $z$, $u$
has a
polynomial growth near
the boundary of the
bidisk $\D\times \D$:
$$
|K_\alpha(\phi,s;  z,u)|
\leqslant \exp\{4\pi\,(1+ |{\rm Im}\, s|\}
(1-|z|)^{-(1+\alpha+2\,|{\rm Re}\,s|)}
(1-|u|)^{-(1+\alpha+2\,|{\rm Re}\,s|)}
$$
Hence,
for fixed $s$, $\phi$, the coefficients
$a_{kl}=a_{kl}(\phi,s)$
of the series
$$K_\alpha(\phi,s;z,u )=\sum a_{kl}(s,\phi)\ov  z^k u^l$$
have the polynomial growth
 as $k,l\to+\infty$.

Indeed, let $q(z,u)$ be a function
in the bidisk satisfying
$$
|q(z,u)|\leqslant C\cdot \delta^{-h}
\qquad \text{for\quad $|z|\le1-\delta$, $|u|\le1-\delta$}
.$$
The Taylor coefficients $b_{kl}$
of $q(z,u)$ are given by the formula
$$b_{kl}=\frac 1{(2\pi i)^2}
\iint\limits_{|z|=1-\delta, |u|=1-\delta}
\frac{q(z)\,dz\,du}{z^{k+1} u^{l+1}}
.$$
Hence,
 $$|b_{kl}|\le {\rm const}\cdot
\delta^{-h} (1-\delta)^{-k-l-2}.
$$
for all $\delta$. We choose
$\delta=h/(h+k+l+2)$ and obtain
polynomial growth for $b_{kl}$.

For the kernel $K_\alpha(z,u)$, we obtain
in this way
the uniform estimates of the form
$|a_{kl}|\leqslant
A\cdot (1+k+l)^\tau$
in each rectangle
\begin{equation}
|{\rm Im}\, s|\leqslant N, \qquad |{\rm Re}\, s|\leqslant M
.\end{equation}

Since the Taylor coefficients $c_{kl}$ for $f\in W$
exponentially decrease,  the series
$$
J_\alpha f(\phi,s)=\sum c_{k,l} a_{k,l}(s,\phi) \cdot
$$
is absolutely convergent
 (see. (4)),
 the summands are holomorphic
 with respect to
$s$, and the series
is uniformly convergent
on rectangles (10).

\smallskip

{\bf 7. Correspondence of differential
operators.}
The operators of the Lie
algebra $sl(2)\oplus sl(2)$
in the space
$H_\alpha\otimes \ov H_\alpha$
have the form
\begin{align*}
L_0^{(z)}&=
z\frac\partial {\partial z}  +\frac\alpha 2 ;&\qquad
L_1^{(z)}&=z^2\frac\partial {\partial z} +\alpha z;&\qquad
L_{-1}^{(z)}&=\frac\partial {\partial z}&\\
L_0^{(u)}&=\ov u\frac\partial {\partial \ov u}  +
\frac\alpha 2 ;&\qquad
L_{1}^{(u)}&=\frac\partial {\partial \ov u};&\qquad
L_{-1}^{(u)}&=\ov u^2\frac\partial {\partial \ov u}
 +\alpha \ov u&
\end{align*}
For us it will be more convenient the following
collection
of the operators
\begin{align*}
L_0&:=L_0^{(z)}-L_0^{(u)};& L_{-1}&:=L_1^{(u)}-L_{-1}^{(z)};
 & L_1:=L_1^{(z)}-L_{-1}^{(u)};
\\
M_0&:=L_0^{(z)}+L_0^{(u)};& M_{-1}&:=L_1^{(u)}+L_{-1}^{(z)};&
   M_1:=L_1^{(z)}+L_{-1}^{(u)};
\end{align*}
The operators $L_0, L_1, L_{-1}$
span the diagonal subalgebra
$sl_2$ in $sl_2\oplus sl_2$.
Their images under the operator $J_\alpha$
are defined by the formulas
\begin{align*}
&
 J_\alpha\circ\Bigl[z \frac\partial {\partial z}-
    \ov u \frac\partial {\partial \ov u} \Bigr] f(\phi,s)=
\bigl[ \frac\partial {i\partial \phi} \circ J_\alpha \bigr]
f(\phi,s)\\
&
J_\alpha\circ\Bigl[ \ov u^2 \frac\partial {\partial\ov u} +\alpha\ov u -
             \frac\partial {\partial z} \Bigr]f(\phi,s)=
  -\Bigl[ e^{i\phi} \frac\partial {i\partial \phi} +
    (\frac 12+is)e^{i\phi}\Bigr] J_\alpha\,\,  f(\phi,s)\\
&
J_\alpha\circ\Bigl[  z^2 \frac\partial {\partial z}  +\alpha z -
                         \frac\partial {\partial \ov u}\Bigr]=
  \Bigl[ e^{-i\phi} \frac\partial {i\partial \phi} -
    (\frac 12+is)e^{-i\phi}\Bigr] J_\alpha\,\,  f(\phi,s)
\end{align*}

These three formulas easily follow
from
 (8).

 \smallskip

{\sc Theorem 2.}
 {\it The unitary operator $J_\alpha$
  transform the operator $M_0$
 to the operator
\begin{multline}
Q_0f(\phi,s)
=-\frac{(-\tfrac12+is)
 (-\alpha+\tfrac12+is)}{2is}f(\phi,s+i)+
\\+
\frac{(\tfrac12+is)(\alpha-\tfrac12+is)}{2is}f(\phi,s-i)
-
\frac{-\alpha+\tfrac12+is}{2is(-\tfrac12+is)}
     \frac{\partial^2}{\partial\phi^2}f(\phi,s+i)
,\end{multline}
i.e., for any $f\in W$ {\rm(}see our
Section {\rm 6 )},
$$Q_0J_\alpha f=J_\alpha M_0f.$$
The operator $M_1$ transforms to
\begin{multline}
Q_1f(\phi,s)=\\
=e^{i\phi}\biggl[\frac{(\tfrac12+is) (-\alpha+\tfrac12+is)}{2is}f(\phi,s+i)
+\frac{(\tfrac12+is)(\alpha-\tfrac12+is)}{2is}f(\phi,s-i)-\\
\frac{-\alpha+\tfrac12+is}{2is(-\tfrac12+is)}\frac{d^2}{d\phi^2}f(\phi,s+i)+
\frac{-\alpha+\tfrac12+is}{-\tfrac12+is}
      \frac{\partial}{i \partial\phi}f(\phi,s+i)\biggr]
,\end{multline}
and the operator $M_{-1}$ transforms to}
\begin{multline}
Q_{-1}f(\phi,s)=\\
=e^{-i\phi}\biggl[-\frac{(\tfrac12+is) (-\alpha+\tfrac12+is)}{2is}f(\phi,s+i)
+\frac{(\tfrac12+is)(\alpha-\tfrac12+is)}{2is}f(\phi,s-i)-\\
\frac{-\alpha+\tfrac12+is}{2is(-\tfrac12+is)}\frac{d^2}{d\phi^2}f(\phi,s+i)
-\frac{-\alpha+\tfrac12+is}{-\tfrac12+is}
      \frac{\partial}{i \partial\phi}f(\phi,s+i)\biggr]
.\end{multline}

{\sc Proof.}
These formulas can be checked by direct calculations.
For instance, let us consider $M_0$.

Obviously, the operator $M_0$
 is selfadjoint in $H_\alpha\otimes\ov H_\alpha$.
Hence
\begin{multline*}
J_\alpha M_0 f(\phi,s)= \\
\iint\limits_{\D\times \D}
  K(\phi,s;\ov z, u)\Bigl[
 \Bigl(z\frac\partial{\partial z}+
 \ov u\frac\partial{\partial \ov u}+\alpha\Bigr)
\Bigr]
  f(z,\ov u)\cdot
   (1-z\ov z)^{\alpha-2} (1-u\ov u)^{\alpha-2}
 \,\{dz\}\,\{du\}=                          \\=
\iint\limits_{\D\times \D}             \Bigl[
 \Bigl(\ov z\frac\partial{\partial\ov z}+
 u\frac\partial{\partial u}+\alpha\Bigr)
\Bigr]
 K(\phi,s;\ov z, u)\cdot
  f(z,\ov u)  (1-z\ov z)^{\alpha-2} (1-u\ov u)^{\alpha-2}
\,\{dz\}\,\{du\}
\end{multline*}
Thus the first statement of the theorem
is equivalent to the identity
\begin{equation}
 \Bigl(u\frac\partial{\partial u}+
\ov z\frac\partial{\partial \ov z}+\alpha\Bigr)
  K(\phi,s;\ov z, u)-
Q_0\bigl[
 K(\phi,s;\ov z, u)\bigr]=0
\end{equation}
After division by $K$, this  identity
transforms to the form
\begin{multline}
K^{-1}
\Bigl\{u\frac\partial{\partial u}+
z{\frac {\partial}{\partial\ov z}\Bigr\} K}
 +\alpha
+K^{-1}\frac{(-\tfrac12+is)
 (-\alpha+\tfrac12+is)}{2is}K(s+i)-
\\-
K^{-1}
\frac{(\tfrac12+is)(\alpha-\tfrac12+is)}{2is}K(s-i)
+
K^{-1}
\frac{-\alpha+\tfrac12+is}{2is(-\tfrac12+is)}
     \frac{\partial^2}{\partial\phi^2}K(s+i)
\end{multline}
The function in the left side
is a long rational expression
 in $z$, $\ov u$, $e^{i\phi}$, $s$;
  the identity can be easily verified by MAPLE.

 Let us explain how to
verify the identity (15) "by hands".
Each summand of  (15)
can be represented as
a linear combination
$$
a(s)+b(s)\cdot \frac 1{1-z\ov u}
+c(s)\cdot\Bigl[\frac 1{1-ze^{i\phi}}+
    \frac1{1-\ov ue^{-i\phi}}\Bigr]+
    d(s)\cdot\frac{(1-ze^{i\phi})(1-\ov u e^{-i\phi})}
    {1-z\ov u}$$
After this, it remains to sum the coeffitients.

\smallskip

{\sc Remark.}
Let us apply the
operator $Q_0$ to
the functions $f(\phi,s)= g(s)$.
Then the equations
$$Q_0 g(s) = (k+\alpha)  g(s)$$
coincide with a
partial case of the difference equations for the continuous
dual Hahn polynomials
 (see \cite{AAR}, (6.10.9)).

\smallskip

{\bf 8. Some remarks.}
In the work of the author \cite{Nerindex},
there were obtained some elements
of an operational calculus for
the index hypergeometric transform
 (it is called also
 by Olevsky transform
or Jacobi transform,
see \cite{Wey}, \cite{Koo})
\begin{equation}
g(x)\mapsto \widehat g(s)=
\frac 1 {\Gamma(b+c)} \int_0^\infty
g(x)\,\,{}_2F_1(b+is,b-is;b+c;-x) x^{b+c-1}(1+x)^{b-c}    dx
,
\end{equation}
In \cite{Nerindex},
it was shown that
the index hypergeometric transform
maps
the differential operators
$$ Ag(x)=x g(x); \qquad Bg(x)=x\frac\partial{\partial x} g(x)$$
(and hence all the
 operators admitting polynomial
expression in
$x$, $x\frac{\partial}{\partial x}$)
to difference operators in imaginary direction;
see also the work of Cherednik \cite{Che}
containing some similar statements
for simmetric functions
in multidimensional case.
Existence of the formulas (11)--(13)
more or less follows
from these results, but this way
for obtaining the expressions (11)--(13)
also is not very simple.

\smallskip

There arises the following question.

\smallskip

{\sc Question.}
Is it possible to write explicitly
operators of the overalgebra
for the case of
$L^2$
on a pseudo-Riemannian symmetric space
and for the kernel representations?

 Is it possible
to do this at least for rank 1
symmetric spaces?

\sf Math.Phys. Group,
Institute of Theretical and Experimental Physics,

B.Cheremushkinskaya, 25, Moscow 117259, Russia

\tt neretin@main.mccme.rssi.ru


\smallskip


\begin{thebibliography}{cc}

\bibitem{AAR}
Andrews, G.R.,  Askey, R., Roy, R.,
{\it Special functions.}, Cambridge Univ. Press, 1999

\bibitem{Ban}
van den Ban E.P., Schlichtkrull H.
{\it  A residue calculus for root systems.}
Compositio Math. 123 (2000), no. 1, 27--72.

\bibitem{Bar}
Bargmann V.
{\it Irreducible unitary representations of
the Lorentz group.} Ann. Math, 48 (1947), 568--640

\bibitem{Ber}
Berezin, F.A. {\it On relations between covariant and contravariant
symbols of operators
for complex classical domains.}
 Dokl. Akad Nauk SSSR, 241, No 1 (1978), 15--17;
English translation: Sov. Math. Dokl. 19 (1978), 786--789

\bibitem{dB}
de Branges, L.
{\it Tensor product spaces.}
J. Math. Anal. Appl. 38 (1972), 109--148.


\bibitem{Che}
Cherednik I.
{\it Harish-Chandra transform and difference operators.}
Preprint, available via http://arXiv.org/abs/math/9706010

\bibitem{Del}
Delorme P. {\it  Formule de Plancherel
 pour les espaces sym\'etriques réductifs},
 Ann. of Math. 147
   (1998), 417--452.

\bibitem{vD}
van Dijk, H., Hille, S.C.,
 {\it Canonical representations related to hyperbolic
spaces},
J.Funct.Anal., 147, 109--139   (1997).






\bibitem{GN} Gelfand I.M., Naimark M.I.,
{\it Unitary representations of
classical groups.}
 Trudy MIAN., t.36 (1950);
 German translation: Gelfand I.N., Neumark M.A.,
 {\it Unitare Darstellungen der klassischen gruppen.},
 Akademie-Verlag, Berlin, 1957.



\bibitem{GK1}          
Gindikin, S.G., Karpelevich, F.I., {\it Plancherel measure for
Riemannian symmetric spaces of non-positive curvature.}
 Dokl. Akad Nauk SSSR,
 145, 252--255(1962); English translation Sov.Mat.Dokl,3, 962--965
 (1962).

\bibitem{GK2}            
Gindikin, S.G., Karpelevich, F.I.,
{\it On an integral connected with symmetric
 Riemannian space of non-positive curvature.}
 Izv. Akad. Nauk SSSR, Ser Mat., 30, 1147--1156(1966)
; English translation in
 Transl.Amer.Math.Soc., 85, 249--258(1969)



\bibitem{Har1}
Harinck P. {\it Plancherel formula pour $\GL(n,\C)/\U(p,q)$.}
J. Reine Angew. Math. 428 (1992), 45--95.

\bibitem{Har2}
Harinck P.
 {\it Fonctions orbitales
  sur $G_\C/G_\R$. Formula d'inversion
des integrales orbitales et formula de Plancherel.}
J. Funct. Anal., 153 (1998), 52--107

\bibitem{HC}
Harish-Chandra,
 {\it Harmonic analysis on real semisimple groups.  III:
The Maas--Selberg relations and the Plancherel formula.}
Ann. Math.,  104 (1976), 117-201.

\bibitem{Hel}           
Helgason, S., {\it Groups and geometric analysis.}
 Acad. Press (1984);

\bibitem{KS}
Koekoek, R., Swarttouw, R.F.,
{\it Askey scheme of hypergeometric
orthogonal polynomials and their $q$-analogues.}
Delft University of Technology, 1994;
available via
{\tt http://aw.twi.tudelft.nl/~koekoek/research.html}.



\bibitem{Koo}
Koornwinder, T.H.,
{\it  Jacobi functions and analysis on noncompact symmetric
spaces.} in {\it Special functions: group theoretical aspects
and applications,}
eds. Askey, R., Koornwinder T.H., Schempp W., 1--85,
D. Reidel Publ. Co., Dordrecht--Boston, 1984.

\bibitem{Mak}
 Makarevich, B. O.,
 {\it Open symmetric orbits
 of reductive groups in symmetric
$R$-spaces.}
 Mat. Sb. (N.S.) 91(133) (1973), 390--401, 472;
 English transl.:
 Math. USSR-Sb., 20 (1973), 406--418.

\bibitem{tensor}
Molchanov V.F.
{\it
Tensor products of unitary representations of the
three-dimensional Lorentz group.}
 Izv. Akad. Nauk SSSR Ser. Mat. 43 (1979), no. 4,
860--891, 967.
English translation: Math USSR Izv., 15 (1980), 113--143.


\bibitem{Molhyper}
Molchanov V.F.
{\it  Plancherel's formula for hyperboloids.}
 Trudy Mat. Inst. Steklov. 147 (1980), 65--85
English translation:
{\it Boundary value problems of mathematical physics. X.}
 Proc. Steklov Inst. Math. 1981, no. 2.
Amer. Math. Soc., Providence, R.I.,
 1981. pp. 63--83.


\bibitem{Molrank}
Molchanov V.F.
{\it The Plancherel formula for
pseudo-Riemannian symmetric spaces of rank $1$.}
Dokl. Akad. Nauk SSSR 290 (1986), no. 3, 545--549.
English translation:
Sov. Math. Dokl. 34 (1987), 323--326

\bibitem{MolVINITI}
Molchanov V.F.
{\it Harmonic analysis on homogeneous spaces}.
 Encyclopaedia Math. Sci., 59,
Representation theory
 and noncommutative harmonic analysis, II, 1--135,
Springer, Berlin, 1995.


\bibitem{pse}
Neretin Yu.A.
{\it Pseudo-Riemannian symmetric spaces:
 single-type realizations and open embeddings into
Grassmannians.}
Zap. Nauchn. Sem. S.-Peterburg.
 Otdel. Mat. Inst. Steklov. (POMI) 256 (1999),
   145--167,
English translation,  J.Math.Sci, New York;
preprint version is available via
 {\tt arXiv.org/abs/math/9905014}





\bibitem{Nerpl}
Neretin, Yu.A.,
{\it Plancherel formula for Berezin deformation
 of $L^2$ on Riemannian
symmetric space}, to appear in J. Funct. Anal.,
 preprint   version is available via
 {\tt  http://arXiv.org/abs/math/9911020}



\bibitem{Nerindex}
Neretin Yu.A.
{\it Index hypergeometric transformation
 and an imitation of the analysis of Berezin
kernels on hyperbolic spaces.}
Mat. Sb. 192 (2001), no. 3, 83--114;
English translation in Sb. Math. 192 (2001),
 no. 3-4, 403--432;
 is available via
 {\tt http://arXiv.org/abs/math/0104035}


\bibitem{Nerhyper}
Neretin Yu.A.
{\it Matrix balls, radial analysis of Berezin
kernels and hypergeometric determinants.}
 Moscow Math.J.,  vol.1, 157 -- 220;
preprint version
 is available via
 {\tt http://arXiv.org/abs/math/0012220}

\bibitem{Puk}
Pukanszky, L., {\it On the Kronecker products of irreducible
 unitary representations  of the
 $2\times2$ real unimodular  group.}
  Trans. Amer. Math. Soc., 100
(1961), 116--152



\bibitem{UU}
Unterberger, A., Upmeier, H., {\it The Berezin transform and
 invariant differential operators}.
Comm.Math.Phys.,164, 563--597(1994)


\bibitem{VGG}
Vershik A.M., Gelfand I.M., Graev M.I.,
{\it Representations of , where $R$ is function ring.}
Uspehi Mat. Nauk 28 (1973), No 5,83--128;
English translation: Russian Math.
Surveys 28 (1973), No 5, 87--132











\bibitem{Wey}
Weyl, H.,
{\it Uber gewonliche lineare Differentialgleichungen
mis singularen Stellen und ihre Eigenfunktionen {\rm 2} Note}.
Nachr. Konig. Gess. Wissen. Gottingen, Math.-Phys.,
1910, 442-467; Reprinted in
Weyl, H., {\it Gessamelte Abhandlungen}, Bd. 1, 224--247,
Springer, 1968.


 \end{thebibliography}
\end{document}